\documentclass[11pt]{article}
\usepackage{amssymb}
\usepackage{amsmath}

\usepackage{amssymb}
\usepackage{amsmath}

\newenvironment{proof}{\noindent {\bf Proof }}
{\hfill $\bullet$ \vspace{0.25cm}}

\def\one{{\bf 1}\hskip-.5mm}

\def\E{{\mathbb E}}
\def\P{{\mathbb P}}
\def\R{{\mathbb R}}
\def\Z{{\mathbb Z}}

\def\ZZ{{\mathbb Z}}

\def\LL{{\mathcal G}}
\def\SS{{ \cal S}}
 \def\F {{\mathcal F}}
\def\FF {{\mathcal F}}

\def\L {{\Lambda}}

\def\ge {{\varepsilon}}
\def\ll{{\bar \Lambda_n}}
\def\gg{{\underline{\gamma}}}

\newtheorem{theo}{Theorem}
\newtheorem{prop}{\indent Proposition}

\newtheorem{rem}{\indent Remark}
\newtheorem{lem}{\indent Lemma}
\newtheorem{defin}{\indent Definition}
\newtheorem{cor}{\indent Corollary}

\title{Partially observed Markov random fields are variable neighborhood random fields }

\author{M. Cassandro \and A.~Galves \and E.~L\"ocherbach}

\date{April 17, 2012}

\begin{document}

\maketitle \begin{abstract}
The present paper has two goals. First to present a natural example of a new class of random fields which are the variable neighborhood random fields. The example we consider is a partially observed nearest neighbor binary Markov random field.  The second goal is to establish sufficient conditions ensuring that the variable neighborhoods are almost surely finite. We discuss the relationship between the almost sure finiteness of the interaction neighborhoods and the presence/absence of phase transition of the underlying Markov random field. In the case where the underlying random field has no phase transition we show that the finiteness of neighborhoods depends on a specific relation between the noise level and the minimum values of the one-point specification of the Markov random field. The case in which there is phase transition is addressed in the frame of the ferromagnetic Ising model. We prove that the existence of infinite interaction neighborhoods depends on the phase.  \end{abstract}



\section{Introduction}

Recent experimental data suggest that populations of neurons have interactions of variable
range. There are reasons to believe that the interaction neighborhood of each neuron is not fixed, but changes as a function of the configuration. Actually the same phenomenon seems to be present at different scales 
when instead of individual neurons local sub-populations of neurons act as interaction unity. 
Justifying the variable neighborhood assumption for fields describing populations
of neurons is an important open question in neuroscience. For a general discussion of the geometry of the neuronal connectivity we refer the reader to Braitenberg and Sch\"utz (1998). Concerning the relationship between the time evolution of the neuronal activity and the 
reaction to external stimulations see MacLean et al. (2005). For a very recent statistical and clinical discussion of the way neighborhood interactions between regions of the brain can change we refer to Wang et al. (2010). Finally, for a mathematical model describing variable range interactions in time rather than in space we refer to Cessac (2011) and the references cited therein.

The above observation suggests to model these kind of interacting systems by a new class
of random fields which are the variable neighborhood random fields. This new class of models is a natural
extension to the case of random fields of the notion of stochastic chains with memory of variable length introduced by Rissanen (1983).

Random fields with variable interaction neighborhoods have recently gained interest, and some papers are devoted to the study of such kind of new models, see Dereudre et al. (2011) and L\"ocherbach and Orlandi (2011). The first paper focusses on the 
problem of existence of these models in $\R^d .$  The second paper addresses the problem of statistical inference, mainly in the case of bounded interaction range. 

The present paper has two goals. First we present a simple and natural class of variable neighborhood random fields, namely the incompletely observed Markov random fields. The second goal is to search for sufficient conditions ensuring that the variable neighborhoods are simultaneously finite for almost every realization of the field. 

The model we consider is a nearest neighbor Markov random field taking the values $+1$ or $-1 .$  At each site there is an independent random mechanism which hides the actual value of the spin and replaces it in the observed data by the value $-1 .$ This can be seen as a black and white picture in which random noise affects the readability of some of the pixels which appear black independently of the actual color. In particular, the noise mechanism introduces a high bias into the system. As in the one-dimensional case where random observations of Markov chains lead to processes having infinite memory, see e.g. Collet and Leonardi (2009), a priori such a model is a random field having infinite memory. However, in this particular case, the partially observed Markov random field is indeed a variable neighborhood random field, and the 
relevant neighborhoods needed in order to determine the spin at a given site will be regions surrounded by a circuit of sites having all spins equal to $+1 .$  This is the content of Theorem \ref{theo:vnrf}. 

Several questions arise naturally in this context. First, is there a relation between presence or absence of phase transition for the underlying random field in $\Z^2 $ and finiteness of the interaction regions of the variable neighborhood random field? Does the absence of phase transition always imply that the interaction regions are finite almost surely? Do infinite interaction regions always exist in the regime of phase transition? It turns out that the question of presence/absence of phase transition and the question of finiteness of interaction neighborhoods are related in a more 
intricate way than we would have guessed naively.

The case in which there is no phase transition is treated in our Theorem \ref{theo:two}. If the minimum values of the one-point specification of the original Markov random field are large enough, we show that
the two situations are possible, depending on the specific relationship between the perturbation level and the specification minima.  The proof relies on a coupling argument. 

The case in which there is phase transition is addressed in the frame of partially observed ferromagnetic Ising models.
We show that, as a consequence of the bias of the noise, the plus phase and the minus phase behave differently when the perturbation level and the temperature are small enough. Namely, in the plus phase all interaction neighborhoods will be finite almost surely, while in the minus phase, infinite interaction regions will always exist with strictly positive probability.      
This is the content of Theorem \ref{theo:main}. The proof of this theorem is based on a Peierls contour counting argument.

This paper is organized as follows. Definitions, notation and main results are presented in Section \ref{section:def}.
The proofs of Theorem \ref{theo:vnrf}, Theorem \ref{theo:two} and Theorem \ref{theo:main}  are presented successively in Sections \ref{section:vnrf}, \ref{section:proofmain} and  \ref{section3}.

\section{Definitions, notation and main results}\label{section:def}
Let $A:= \{ -1 ,1\}$ and $S =A^{\Z^2} $ be the set of all possible
configurations. We endow $S$ with the product sigma algebra ${\cal S} .$  
Fixed configurations will be denoted by lowercase letters $x, y, z .$ 
A point $i \in \Z^2 $ is called a site. 

If $ x \in S$ is a configuration, then for any $i \in \Z^2 ,$ $x_i$ will denote the value of the configuration at site $i.$ Given
a subset $F \subset \Z^2,$ we will also denote $x_F = \{ x_i , i \in F \}. $ 
Let $X = \{ X_i :  i \in \Z^2\}  $ be the collection of projections on $S ,$ defined by $ X_i (x) = x_i $
for all $i \in \Z^2 .$ We introduce the  following $\sigma-$algebras:  For any $\Lambda \subset \Z^2, $ 
let 
$$ {\cal F}_\Lambda = \sigma \{ X_i : i \in \Lambda \}  .$$ 

\begin{defin}
Any  probability measure on $(S, {\cal S})$ will be called a random field.
\end{defin}
We consider random fields on $(S, {\cal S})$ which are defined by their local specifications, see Dobrushin (1970). In order to do so, we
recall the notion of specification from Georgii (1988). 

\begin{defin}  A  specification on $(S, {\cal S})$  is a family $ P= \{P_{\Lambda}\}_{\Lambda \subset  \Z^2}  $ of probability kernels on $(S, {\cal S})$ such that  
  \begin{itemize}
  \item[(a)] For each $\Lambda \subset \Z^2$ finite and each $ B \in {\cal S}$, the
    function $ P_{\Lambda}(B \mid\cdot\,)$ is $
    \FF_{\Lambda^c}-$measurable.

\item[(b)] For each $\Lambda\subset \Z^2$ finite and each $ B \in
  \FF_{\Lambda^c}$, $ P_{\Lambda}(B\mid y)= 1_{B}(y ).$

\item[(c)] For any pair of finite subsets $ \Lambda$ and $\Delta$, with
  $\Lambda \subset \Delta\subset\Z^2$, and any measurable set $B$,
 \begin{equation}
 \int P_\Lambda(B\mid z ) \, P_{\Delta}(d z\mid  y) 
 \;=\; P_{\Delta}(B\mid y)
 \label{eq:3}
\end{equation}
for all $y \in S$.
\end{itemize}
\end{defin} 
From $(a)$ and $ (b)$ above it follows that $ P_\Lambda $ can be identified with probability weights $ p_\Lambda ( x_\Lambda  | y_{\Lambda^c } ) $ such that for any $ B \in {\cal F}_\Lambda ,$ 
$$ P_\Lambda (B | y) = \sum_{ x_\Lambda : (x_\Lambda, y_{ \Lambda^c})  \in B} p_\Lambda ( x_\Lambda  | y_{\Lambda^c } ) .$$

\noindent
\begin{defin}
 A random field $\mu$ on $(S,\SS)$ is {consistent} with
    a specification $P $ if for each finite subset $\Lambda\subset\Z^2  ,$
    \begin{equation}
  \int \mu(dx) \; P_\Lambda(B \mid x)
  \;=\; \mu(B) ,
  \label{eq:3.1}
  \end{equation}
  for every $B\in \SS $.
We write ${\cal G} (P)$ for the set of all random fields consistent with the specification $P .$
\end{defin}
In the sequel, if $P$ is a specification, for $ \Lambda = \{ i \} , $ instead of writing $p_{\{ i \}} ( \cdot | \cdot )  $ for the probability weights, we shall use the
short-hand notation $ p_i ( \cdot | \cdot ) .$ 

One important class of random fields are the Markov fields. For the reader's convenience we recall here the basic definitions. 
The $L^1 -$norm $ \| \cdot \|_ 1 $ 
is defined as follows.  
If $ i = ( i_1, i_2) \in \Z^2  ,$ then 
$$ \| i \|_1 = | i_1 | + | i_2| . 
$$ 
For any subset $F \subset \Z^2, $ the boundary $ \partial F $ is defined as 
$$ \partial F = \{ j \in F : \exists k \in \Z^2 \setminus F : \| j - k \|_1 = 1 \} .$$

\begin{defin}
Two points $ i $ and $ j \in \Z^2 $ are called $ L^1 -$neighbors 
if $ \| i - j \|_1 = 1 .$ 
\end{defin}

\begin{defin}
Let $P$ be a specification. $\mu \in {\cal G} (P)  $ is
a Markov random  field of order $1$ if for all $ i \in \Z^2 ,$ the function $  P_{i  } (\{ X_i = + 1 \} | \cdot  ) $ is ${\cal F}_{ \partial (\{ i \}^c) }-$measurable, where $ \partial (\{ i \}^c) = \{ j \in \Z^2  : \| i - j \|_1 = 1 \} .$ 
\end{defin}

We now introduce the notion of variable neighborhood random fields.  
\begin{defin}
Let $P$ be a specification and $ \P \in {\cal G} (P) .$ We say that $ \P$ is a 
variable neighborhood random field if for any $ i \in \Z^2 $ there exists a mapping $C_i : A^{\Z^2 \setminus \{i \} } \to {\cal P} ( \Z^2 \setminus \{ i \} ) $ such that the following statements hold.
\begin{enumerate}
\item
For any site $ i \in \Z^2 $ and any subset $ \Lambda \subset \Z^2 ,$ the event $ \{ C_i (X) = \Lambda \} $ belongs to $ { \cal F}_\Lambda .$ 
\item
For all $i \in \Z^2 ,$  the map $x \mapsto   P_{ i } ( \{ + 1 \} | x )  $ is $ { \cal F}_{C_i (X) }-$measurable.  Here, 
$$  { \cal F}_{C_i (X) } = \{ B \in {\cal S} :   B \cap \{ C_i (X) = \Lambda \} \in {\cal F}_\Lambda ,\mbox{ for all } \Lambda \subset \Z^2\} .$$
\item
$C_i (x)$ cannot be shortened. This means that if there is another collection of maps $\tilde C_i , i \in \Z^2 ,  $ such that the above properties hold, then $ C_i (x) \subset \tilde C_i (x ) $ for all $ i $ and $x.$ 
\end{enumerate} 
We call $ x_{C_i(x)} $ the context of site $i,$ given the configuration $x.$ The functions $\{C_i (\cdot), i \in \Z^2 \} $ will be called context support functions of the variable neighborhood random field. 
\end{defin}

From now on we shall write shortly VNRF for variable neighborhood random field.

The goal of this paper is to introduce a natural example of VNRF's, namely incompletely observed Markov random fields.
By this we mean the following. For any fixed $ \varepsilon \in ] 0, 1[ ,$ let $ \nu_\varepsilon$ be the Bernoulli product measure 
\begin{equation}\label{eq:bernoulli}
\nu_\varepsilon = \prod_{i \in \Z^2} ( \varepsilon \delta_{-1} + (1 - \varepsilon ) \delta_{+1} ) 
\end{equation}
on $(S, {\cal S} ) .$ Thus under $ \nu_\varepsilon,$ the coordinates $ X_i, $ $i \in \Z^2 ,$ are i.i.d. random variables taking the value $+ 1$ with probability $ 1 - \varepsilon $ and the value $-1$ with probability $\varepsilon.$ For each site $i, $ its original color chosen according to $\mu $ will be observed only with probability $ 1 - \ge , $ and with probability $\ge ,$ we loose any information concerning the color and report as output the value $-1.$ 

Mathematically speaking, this means the following. For any measure $\mu \in {\cal G} (P)$ we consider the product measure $ \mu \otimes \nu_\varepsilon $ on $ ( S \times S , {\cal S} \otimes {\cal S})$ and consider the probability measure $ \P^{ \varepsilon } $ on $ ( S , {\cal S}) $ which is the image measure of $ \mu \otimes \nu_\varepsilon $ under the operation of taking the point-wise minimum
\begin{equation}\label{eq:pepsilon}
 S \times S \ni (x^1 , x^2) \mapsto x \in S  : \mbox{ for all } i \in \Z^2, \; x_i = x^1_i  \wedge x^2_i .
\end{equation}
In other words, observing a realization of the random field $ \P^{  \varepsilon }  $ amounts to saying that we observe a realization of the original random field $\mu,$ where for each site, independently of the original value of the field and independently of the other sites, its value is replaced by the value $-1 .$  

It turns out that incompletely observed random fields as defined above are VNRF's. This is the content of our first theorem. Before stating it, we recall some definitions. For any finite subset $F \subset \Z^2, $ the interior $\mbox{\it \r{ F}} $ is given by 
$$ \mbox{\it \r{ F}} = F \setminus \partial F .$$

\begin{defin}
A self-avoiding path in $\Z^2 $ is a finite sequence $ \Gamma = (i_1 , \ldots , i_n ) $ of sites such that $ i_j $ and $ i_k $ are $L^1 -$neighbors if and only if $ | j - k | = 1 .$
\end{defin}

Recall also that a set $ F \subset \Z^2$ is called $L^1-$connected if for any pair of points $ i , i'  \in F ,$ $i \neq i',$ there exists an integer $ n \geq 1 $ and a self-avoiding path $(i_1 , \ldots , i_n ) $ of sites in $F$ such that $i_1 = i $ and $i_n = i' .$   

%

The following theorem holds. 
     
\begin{theo}\label{theo:vnrf}
If $\mu \in {\cal G} (P) $ is a Markov random field of order $1$ and $\nu_\varepsilon$ the Bernoulli probability measure of (\ref{eq:bernoulli}), then the measure $\P^\ge $ defined in (\ref{eq:pepsilon}) is a VNRF, and the context support functions are given by  
$$ C_i (x) = \left( \cap \{ F \subset \Z^2 :   i \in  \mbox{\it \r{ F}}  , \; F  \mbox{ is  $L^1-$connected and   }    x_{ \partial  (F)} \equiv + 1 \} \right) \setminus \{ i \} ,$$
if there exists at least a finite set satisfying these conditions. Otherwise, we put $ C_i (x) = \Z^2 \setminus \{ i \} .$ 
\end{theo}

It is natural to ask whether for a given Markov random field model $\mu $ and a given $\ge ,$ all contexts $C_i (x) $ will be finite almost surely or not. In the regime of absence of phase transition, a first answer can be given by using coupling arguments. Call 
\begin{equation}\label{eq:lambdaplus}
\lambda_0^+ = \inf_{i \in \Z^2} \inf_{ x \in S } p_i ( +1 | x ) ,  \quad \lambda_0^- = \inf_{i \in \Z^2} \inf_{ x \in S } p_i ( -1 | x ) .\end{equation}

Let $p^*$ be the critical probability for the site percolation model in $\Z^2 .$ For a general presentation of percolation models we refer the reader to the classical treatise by Grimmett (1999). Then we have the following result. 

\begin{theo}\label{theo:two}
Grant the conditions of Theorem \ref{theo:vnrf}. If 
\begin{equation}\label{eq:conditionbetapetit}
( 1 - \ge ) \lambda_0^+ > 1 - p^* , 
\end{equation}
then
$$ \P^{\ge} \left( \bigcap_{i \in \Z^2} \left\{ |C_i (X)|  < \infty \right\}   \right) = 1  .$$
On the other hand, if 
\begin{equation}\label{eq:conditionbetapetit2}
\ge + ( 1 - \ge) \lambda_0^- > p^* ,
\end{equation}
then
$$  \P^{ \ge }  \left(\bigcup_{i \in \Z^2 } \left\{  | C_i (X)| = \infty \right\}  \right) > 0 .$$ 
\end{theo}
In the above statements, $| C_i (X) | $ means the cardinal of the set $ C_i (X) .$ We give the proof of this theorem in Section \ref{section3} below.

We now consider the regime in which there is phase transition. We address the question of finiteness of contexts in the framework of the ferromagnetic Ising model. 
\begin{defin}\label{ex:exa1}
The homogeneous ferromagnetic Ising model is defined by the following specification.
For any $\beta \geq 0 ,$ $x, y \in S$ and any finite subset $ \Lambda \subset \Z^2,$ 
\begin{equation}\label{eq:Ising}
 p_{\Lambda , \beta} ( x_\Lambda | y_{\Lambda^c}  ) = \frac{ 1 } { Z^{y}_{\Lambda , \beta} } \exp \left(   \beta \left[  \, \frac12 \sum \limits_{i \in \Lambda } \sum \limits_{j \in \Lambda , \| i - j \|_1 = 1  } x_i x_j + \sum \limits_{i \in \Lambda} \sum \limits_{ j \in \Lambda^c ,\| i - j \|_1 = 1} x_i y_j   \right]
 \right)  , 
\end{equation}
where 
$$  Z^{y}_{\Lambda , \beta } = \sum \limits_{ x_\Lambda \in A^{\Lambda} } \exp \left(   \beta \left[  \frac12 \sum \limits_{i \in \Lambda }  \sum \limits_{j \in \Lambda , \| i - j \|_1 = 1  }x_i x_j + \sum \limits_{i \in \Lambda} \sum \limits_{ j \in \Lambda^c , \| i - j \|_1 = 1} x_i y_j \right] \right) .$$
\end{defin}
It is well known, see for instance Georgii (1988) and Presutti (2009), that that there exists a critical value $\beta^c ,$ such that for all $\beta > \beta^c $ the set ${\cal G} (P)  $ contains two extremal measures 
$\mu_\beta^-$ and $\mu_\beta^+ $ which are the pure states obtained by passing to the limit $\Lambda \to \Z^2 ,$ taking the external configuration $ y_j = -1 $ for all $ j \in \Lambda^c $ and  $ y_j = +1 $ for all $ j \in \Lambda^c,$ respectively. 

We write $ \P^+_{  \beta, \ge} $ for the image measure of $ \mu^+_\beta  \otimes \nu_\varepsilon $ under the map
$$ S \times S \ni (x^1 , x^2) \mapsto x \in S  : \mbox{ for all } i \in \Z^2, \; x_i = x^1_i  \wedge x^2_i .$$ 
$\P^-_{ \beta, \ge } $ 
is defined in an analogous way. 

We have the following result.
\begin{theo}\label{theo:main}
The following two statements hold. \\
1. 
For sufficiently large $ \beta > \frac12 \, \ln 3  $ and $\ge <\frac13 - e^{ - 2 \beta}  ,$ 
$$ \P^{+}_{ \beta , \ge } \left( \bigcap_{i \in \Z^2} \left\{ |C_i (X)|  < \infty \right\}   \right) = 1  .$$
2. For all $ \beta > \beta^c ,$ 
$$  \P^{ -}_{ \beta , \ge } \left(\bigcup_{i \in \Z^2 } \left\{  | C_i (X)| = \infty \right\}  \right) = 1  .$$ 
\end{theo}

\begin{rem}
The high-temperature results of Theorem
 \ref{theo:two} apply very nicely in the framework of the ferromagnetic Ising model. In this case, 
$$ \lambda_0^+ = \lambda_0^- = ( 1 + \exp ( 8 \beta ) )^{-1 },$$
and for example condition (\ref{eq:conditionbetapetit}) reads as 
$$ \ge < p^*  \mbox{   and    } \beta < \frac18 \ln \left( \frac{1 - \ge}{ 1 - p^*} - 1 \right).$$
\end{rem}

\section{Proof of Theorem \ref{theo:vnrf}}\label{section:vnrf}

Write $ P^\ge $ for the specification of $\P^\ge $ and $p^\ge_i$ for the associated one-point conditional probabilities. The strategy of our proof is to provide an explicit formula for $ p^\ge_i ( +1 | x ) $ (see (\ref{eq:version}) below) which shows that $ p^\ge_i ( +1 | x ) $ depends only on $x_{C_i(x)} .$  

Let $\Lambda_n = [ - n, n ]^2 .$ Write $ \bar C_i ( X) = C_i (X) \cup \{ i \} .$ Let $ B  \in {\cal S}.$ We only have to consider the event $ \{ | C_i (X) | < \infty \} . $ We start by evaluating
$$
 \P^\ge (  B ;  | C_i (X)| < \infty ) 
=     \lim_n \P^\ge ( B;   C_i (X) \subset \Lambda_n ) .
$$
We have for any fixed $n,$
$$ \P^\ge ( B ; C_i (X) \subset \Lambda_n ) = \sum_{ F \subset \Lambda_n } 
\P^\ge ( B  ; \bar C_i (X) = F   ).$$
Here, we sum over all sets $F \subset \Lambda_n$ which are $L^1-$connected and contain $i$ in their interior. We evaluate each of these terms $ \P^\ge (  B ;   \bar C_i (X) = F ) ,$ for a fixed set $F .$ In order to simplify notation, we write $  \Lambda = \mbox{\it \r{ F}} $ for the interior of $F.$  Notice that $ i \in \Lambda .$ 

Observe that 
$$ \P^\ge (  B ;   \bar C_i (X) = F )  = \sum_{ z_F \in A^F  : \bar C_i (z) = F } \P^\ge ( B , X_F = z_F ) .$$
For any fixed $F$ and $z_F, $ notice that we can rewrite $ 1_B(x) 1_{z_F}(x_F) = 1_{ B^{z_F}}(x) 1_{ z_F }(x_F),$ where 
$$ B^{z_F} = \{ x \in S : (z_F, x_{F^c} ) \in B \} \in {\cal F }_{F^c } .$$

Hence, using Fubini's theorem and since $ z_{\partial F} \equiv + 1 , $ by construction of $ \P^\ge ,$ 
\begin{multline*}
\P^\ge (B;  X_F = z_F )
  = \int_S \nu_\ge (dx^1) 1_{\{ x^1_{\partial F} \equiv +1 \} }
\\
 \left[ \int_S \mu (dx^2 ) 1_{ \{x^2_{\partial F}  \equiv +1\}} 
  1_{  B^{z_F} } ( x^1 \wedge x^2)  1_{ z_\Lambda } ( x^1 \wedge x^2 )_\Lambda  \right]  ,
\end{multline*}
where $ x^1 \wedge x^2 $ denotes the configuration
$$( x^1 \wedge x^2 )(i) = x^1 (i ) \wedge x^2 (i) , \mbox{ for all } i \in \Z^2 .$$

Observe that for fixed $x^1 ,$ the mapping $ x^2 \mapsto 1_{  B^{z_F} } ( x^1 \wedge x^2 )$ is ${\cal F}_{ F^c}-$measurable.
Hence we obtain, applying (\ref{eq:3.1}), for a fixed configuration $x^1 ,$ 
\begin{multline*}
 \int_S \mu (dx^2 ) 1_{ \{x^2_{\partial F} \equiv +1\}} 
  1_{ B^{z_F} } ( x^1 \wedge x^2) 1_{ z_\Lambda } ( x^1 \wedge x^2 )_\Lambda    \\
= \int_S \mu ( dy ) 1_{ \{y_{\partial F} \equiv +1\}}   1_{  B^{z_F}} ( x^1 \wedge y)   \int_{A^{\Lambda} } p_\Lambda  (d u_\Lambda |y) 1_{ z_\Lambda } ( x^1 \wedge u )_\Lambda  \\
= \int_S \mu ( dy ) 1_{ \{y_{\partial F} \equiv +1\}}   1_{  B^{z_F} } ( x^1 \wedge y)   \int_{A^{\Lambda} } p_\Lambda  (d u_\Lambda |y_{\partial F}) 1_{ z_\Lambda } ( x^1 \wedge u )_\Lambda   \\
=  \int_S \mu ( dy ) 1_{ \{y_{\partial F} \equiv +1\}}   1_{  B^{z_F} } ( x^1 \wedge y)   \int_{A^{\Lambda} } p_\Lambda  (d u_\Lambda |+) 1_{ z_\Lambda } ( x^1 \wedge u )_\Lambda   ,
\end{multline*}
where we used that $\mu$ is a Markov random field of order $1$ and the fact that $ \{ j \in \Z^2 : dist ( j , \Lambda ) = 1 \} = {\partial F} ,$ where $dist $ is the distance defined by the $L^1 -$norm on $\Z^2 .$ 
Putting things together, we obtain
\begin{multline*}
\P^\ge (  B ; X_F = z_F  ) \\
= \int_{S \times S} \nu_\ge (dy^1) \mu (d y^2) 1_{\{ y^1_{\partial F} \equiv y^2_{\partial F} \equiv +1 \}} 1_{ B^{z_F}} ( y^1 \wedge y^2 ) \cdot \\
\left[ \int_{A^\Lambda} \prod_{ j \in \Lambda } [ \ge \delta_{ - 1 } + ( 1 - \ge) \delta_{+1} ] (dx^1_j) 
 \int_{A^\Lambda } p_\Lambda  (d x^2_\Lambda |+ ) 1_{  z_\Lambda} ( x^1 \wedge x^2 )   \right] .
\end{multline*}
Recall that $ z_{\partial F} \equiv + 1 .$ Write for simplicity 
$$
 \varphi_i ( \Lambda , z_F  ) =
 \int_{A^\Lambda} \prod_{ j \in \Lambda } [ \ge \delta_{ - 1 } + ( 1 - \ge) \delta_{+1} ] (dx^1_j) 
 \int_{A^\Lambda } p_\Lambda  (d x^2_\Lambda |+ ) 1_{  z_\Lambda} ( x^1 \wedge x^2 )  .
$$
Thus, 
\begin{multline*}
 \P^\ge ( B ; X_F = z_F ) \\
=  \int_{S \times S} \nu_\ge (dy^1) \mu (d y^2) 1_{\{ y^1_{\partial F} \equiv y^2_{\partial F} \equiv +1 \}}  1_{  B^{z_F} } ( y^1 \wedge y^2 ) \; \varphi_i ( \Lambda , z_F  ) \\
= \int_S \P^\ge ( dy) 1_{\{ y_{\partial F} \equiv + 1 \} } 1_{  B^{z_F} } ( y) \varphi_i ( \Lambda , z_F ) .
\end{multline*}

Therefore, if we define 
\begin{equation}\label{eq:version}
p_i^\ge ( + 1 | z ) := \frac{\varphi_i ( \Lambda, ( +1 , z_{ F \setminus \{ i \} }  ))  } {  \varphi_i ( \Lambda, ( +1 , z_{ F \setminus \{ i \} }  )) + 
\varphi_i ( \Lambda, ( -1 , z_{ F \setminus \{ i \} }  ))  \} )}   ,
\end{equation}
on $\{ \bar C_i (z) = F  \} , $ then we have, applying the above arguments to $ B \cap \{ X_i = +1 \}  $ with $B \in {\cal F}_{ \{ i \}^c } ,$ that 
$$
\int_S \P^\ge ( d x) 1_{ B } ( x)  1_{  z_{F  \setminus \{ i \} }} ( x_{F \setminus \{i\} } )  1_{ x_i = + 1 } 
= \int_S \P^\ge ( d x) 1_{ B } ( x)  1_{  z_{F  \setminus \{ i \} }} ( x_{F \setminus \{i\} } )  \;  p_i^\ge ( + 1 | z ) .
$$
Hence the object defined in (\ref{eq:version}) is a version of the conditional probability $ \P^\ge ( X_ i = +1 | z) $ on the event that $ \bar C_i (z) = F .$ It is evident that $ z \mapsto p_i^\ge ( + 1 | z )  $ is $ { \cal F}_{ C_i (X) }- $measurable. This concludes the proof.

\section{Proof of Theorem \ref{theo:main}}\label{section:proofmain}
This section is devoted to the proof of Theorem \ref{theo:main}. The following lemma is the key of our proof.

\begin{lem}\label{prop:1}
For sufficiently large $ \beta  $ and for all $\ge , $ 
$$ \P^+_{\beta, \ge} \left( \bigcap_{ i \in \Gamma} \{ X_i = - 1 \} \right) \le    \left(  e^{ - 2 \beta} + \ge  \right)^{ |\Gamma| } $$
holds for any path $ \Gamma \subset \Z^2 .$
\end{lem}

In order to give the proof of the lemma, we need to recall the notion of contour. We first recall the notion of dual.

\begin{defin}
We call dual of $ \Z^2 $ the set $(\Z^2) ' =  \Z^2 + ( \frac12 , \frac12 ) .$ For any finite set $ \Lambda \subset \Z^2, $ the dual $ \Lambda ' $ of $\Lambda $ is given by 
$$ \Lambda ' = \{ r = ( r_1, r_2)  \in (\Z^2)'  :  \exists i = (i_1, i_2) \in \Lambda , r_1 = i_1 \pm \frac12, r_2 = i_2 \pm \frac12 \} .$$ 
\end{defin}

We define the $L^1-$norm and the notion of $L^1-$neighbors on the dual in exactly the same way as we did for $\Z^2 .$

A contour is defined as follows. 
\begin{defin}
Any finite sequence $ \gamma = ( r_1 , \ldots , r_n) $ of points in the dual $(\Z^2)' $ is called a contour if $(r_j, \ldots , r_n, r_1, \ldots , r_{ j- 2} ) $ is a path for all $ j = 1 , \ldots , n .$ We write $  | \gamma | = n $ for the length of the contour. We say that two contours $\gamma$ and $\gamma'$ do not intersect if and only if either $R (\gamma )  \cap R (\gamma') = \emptyset  $ or $ \gamma \cap \gamma ' = \emptyset .$
\end{defin}

Joining two neighboring points of the contour $\gamma $ with a straight line, we obtain a closed finite curve in $\R^2 .$ We denote the intersection of its interior with $\Z^2 $ by $ R( \gamma ) .$ 

We have now the elements to prove Lemma \ref{prop:1}.

\begin{proof}{\bf of Lemma \ref{prop:1}.}\\
To any $ (x^1, x^2 ) \in S^2 ,$ we associate the configuration $x = x^1 \wedge x^2 .$  Then
$$ 1_{ \{ x_i = - 1 \} } = 1 - 1_{\{ x_i^1 = +1 \}} + ( 1 - 1_{\{ x_i^2 = +1 \}} ) 1_{\{ x_i^1 = + 1 \}} .  $$
By definition of $ \P^+_{ \beta , \ge } ,$ we have 
\begin{multline}\label{eq:o}
 \P^+_{\beta, \ge} \left( \bigcap_{ i \in \Gamma} \{ X_i = - 1 \} \right)  = 
\mu_\beta^+ \otimes \nu_\ge \left(  \prod_{ i \in \Gamma } \left[ 1 - 1_{\{ x_i^1 = +1 \}} + ( 1 - 1_{\{ x_i^2 = +1 \}} ) 1_{\{ x_i^1 = + 1 \}} \right] \right) \\
= \mu_\beta^+  \left(  \prod_{ i \in \Gamma } \left[1 - ( 1 - \ge)  1_{\{ x_i= +1 \}}    \right] \right) \\
= \sum_{ C \subset \Gamma } \mu_\beta^+ \left( x_C = - 1 , x_{\Gamma \setminus C } = + 1 \right) \ge^{|\Gamma | - | C|} .
\end{multline}

We use that
$$ 
\mu_\beta^+ = \lim_{ \Lambda \to \Z^2 } \mu_{ \beta, \Lambda }^+ , 
$$ where $ \mu_{\beta , \Lambda }^+ = P_{\beta , \Lambda} ( \cdot | +_{\Lambda^c } )$, see (\ref{eq:Ising}),  and where $ +_{\L^c}$ denotes the configuration $ y_j = + 1 $ for all $ j \in \Lambda^c .$ In the sequel we will study the properties of the finite volume measure $ \mu_{\beta, \Lambda}^+ $ and get estimates uniform in $\Lambda ,$ for all $\L $ containing $\Gamma. $

We have
\begin{equation}\label{eq:one}
 \mu_{\beta, \Lambda}^+ \left( x_C = - 1 , x_{\Gamma \setminus C } = + 1 \right)
 = \frac{1}{Z_{\beta , \Lambda}^+}  \sum_{\underline{\gamma}\in \Omega_\L} \prod_{ \gamma \in \underline{\gamma}} e^{ - 2 \beta | \gamma |}  1_{\{ x_C (\underline{\gamma } ) =  - 1\} } 1_{\{ x_{\Gamma \setminus C } (\underline{\gamma }) = + 1 \}} ,
\end{equation}
where
$$\Omega_\L  = \{ \underline{\gamma } = \{ \gamma_1 , \ldots , \gamma_n \} , n \geq 1 , \gamma_i \in \L' \;\mbox{ for all } i , \gamma_1 , \ldots , \gamma_n \mbox{ non intersecting} \} ,$$
and where for any given set of contours $ \underline{\gamma} \in \Omega_\Lambda ,$ $x( \gg) \in A^\Lambda $ denotes the associated configuration. In the above formula we used the classical correspondance between configurations and sets of non intersecting contours, see e.g. Presutti (2009).

In order to evaluate (\ref{eq:one}), let $ C = C_1 \cup \ldots \cup C_n $ be the decomposition of $C$ into the union of its connected components. This means that each $C_i$ is a $L^1-$connected set and $ dist (C_i, C_j ) \geq 2 $ for all $ i \neq j . $ All components $C_1, \ldots , C_n$ must be contained in a contour. Observe that one such contour can contain several components. 
More precisely, for any $m $ between $1 $ and $n,$ let $\{ J_j , j = 1 , \ldots  , m \} $ be a (disjoint) partition of $\{ 1 , \ldots , n \} $
and let 
$$ P_j =  \bigcup_{i \in J_j} C_i .$$ 
Each of the $P_j, 1 \le j \le m ,$ will be surrounded by exactly one contour $ \gamma_j .$
The contours surrounding different $P_j$'s have to be non-intersecting. Moreover, the contours $\gamma_1 , \ldots , \gamma_m$ have to be the only contours that intersect the path $\Gamma .$ 

Given the contours $\gamma_1 , \ldots , \gamma_m$ as in the last paragraph, write
\begin{multline*}
 \Omega ( \L \setminus  \gamma_1 \cup \ldots \cup \gamma_m )  = 
 \{ \underline{\gamma '} = \{ \gamma_1 ' , \ldots , \gamma_k ' \} \in \Omega_\L  : \mbox{ for all } 1 \le i \le k , 1 \le j \le m ,  \\ 
  \gamma_i '  \cap  {\gamma_j }  = \emptyset \mbox{ and } R ( \gamma_i ') \cap \Gamma = \emptyset  \} ,
\end{multline*} 
for the set of all contours not intersecting with $\Gamma$ nor with any of the $ {\gamma_j }, 1 \le j \le m .$ Now we can rewrite (\ref{eq:one}) as follows. 
\begin{eqnarray*}
&&  \frac{1}{Z_{\beta , \Lambda}^+}  \sum_{\underline{\gamma}} \prod_{ \gamma \in \underline{\gamma}} e^{ - 2 \beta | \gamma |}  1_{\{ x_C (\underline{\gamma } ) =  - 1\} } 1_{\{ x_{\Gamma \setminus C } (\underline{\gamma }) = + 1 \}}\\
&&=  \frac{1}{Z_{\beta , \Lambda}^+} \sum_{ m = 1 }^n \sum_{ P_1 , \ldots , P_m} \sum_{ \gamma_1 : P_1 \subset R ( \gamma_1)} e^{ - 2 \beta | \gamma_1|} \ldots \sum_{ \gamma_m : P_m \subset R ( \gamma_m ) } e^{ - 2 \beta | \gamma_m|} 1_{ \{ \gamma_1, \ldots , \gamma_m \mbox{ \tiny non-intersecting } \}} \\
&&  \quad  \quad \quad \quad \quad \quad  \quad \quad \quad \quad \quad \quad \quad \quad \quad \quad \quad \sum_{ \underline{\gamma } \in \Omega ( \L \setminus  \gamma_1 \cup \ldots \cup \gamma_m)  }  \prod_{ \gamma \in \underline{\gamma}} e^{ - 2 \beta |\gamma|} .
\end{eqnarray*} 
Since for any fixed set of $\gamma_1, \ldots , \gamma_m,$  
$$ Z_{\beta , \Lambda}^+ \geq \sum_{ \underline{\gamma } \in \Omega ( \L \setminus  \gamma_1 \cup \ldots \cup \gamma_m)  }  \prod_{ \gamma \in \underline{\gamma}} e^{ - 2 \beta |\gamma|} ,$$
we have
\begin{multline*}
 \mu_{\beta, \Lambda}^+ \left( x_C = - 1 , x_{\Gamma \setminus C } = + 1 \right) \\
\le \sum_{ m = 1 }^n \sum_{ P_1 , \ldots , P_m}  \sum_{ \gamma_1 : P_1 \subset R ( \gamma_1)} e^{ - 2 \beta | \gamma_1|} \ldots \sum_{ \gamma_m : P_m \subset R ( \gamma_m )} e^{ - 2 \beta | \gamma_m|} 1_{ \{ \gamma_1, \ldots , \gamma_m \mbox{ \tiny non-intersecting } \}}.
\end{multline*}

Observe that  
$$\sum_{ i= 1}^m  \min_{ \gamma_i : P_i \subset R ( \gamma_i ) } | \gamma_i | \geq  2  \sum_{i= 1}^m | P_i | 
 +2n = 2 | C| + 2 n\,. $$
As a consequence,
\begin{multline*}
\sum_{ P_1 , \ldots , P_m}   \sum_{ \gamma_1 : P_1 \subset R ( \gamma_1)} e^{ - 2 \beta | \gamma_1|} \ldots \sum_{ \gamma_m : P_m \subset R ( \gamma_m ) } e^{ - 2 \beta | \gamma_m|} 1_{ \{ \gamma_1, \ldots , \gamma_m \mbox{ \tiny non-intersecting } \}} \\
\le 
e^{ - 2 \beta n } e^{ - 2 \beta | C | } \sum_{ P_1 , \ldots , P_m}  \sum_{ \gamma_1 : P_1 \subset R ( \gamma_1)} e^{ -  \beta | \gamma_1|} \ldots \sum_{ \gamma_m : P_m \subset R ( \gamma_m ) } e^{ -  \beta | \gamma_m|}1_{ \{ \gamma_1, \ldots , \gamma_m \mbox{ \tiny non-intersecting } \}} .
\end{multline*}
To obtain an upper bound of the right hand side of the above inequality, we use that
\begin{multline}
 \bigcup_{ P_1, \ldots , P_m } \Big\{ \{ \gamma_1, \ldots , \gamma_m \} : P_1 \subset R ( \gamma_1), \ldots , P_m \subset R ( \gamma_m ) , \gamma_1 , \ldots , \gamma_m \mbox{ non-intersecting} \Big\} \\
= \Big\{ \{ \gamma_1 , \ldots , \gamma_m\} : l_1 \in R ( \gamma_1),    \min \{ l_i : l_i \notin R (\gamma_1) \} \in R ( \gamma _2 ) , \ldots ,  \\
  \min \{ l_i : l_i \notin R ( \gamma_1 ) \cup \ldots \cup R ( \gamma_{ m-1 } ) \} \in R ( \gamma_m), \\
\gamma_1 , \ldots , \gamma_m \mbox{ non-intersecting and for all }  \; 1 \le i \le n , \exists j : C_i \in R ( \gamma_j)  \Big\} ,
\end{multline}
where for every $i = 1, \ldots , n ,$ $ l_i$ is a fixed but otherwise arbitrary element of $C_i .$
Hence
\begin{multline}\label{eq:lastdisplay}
\sum_{ P_1 , \ldots , P_m}   \sum_{ \gamma_1 : P_1 \subset R ( \gamma_1)} e^{ -  \beta | \gamma_1|} \ldots \sum_{ \gamma_m : P_m \subset R ( \gamma_m )} e^{ -  \beta | \gamma_m|} 1_{ \{ \gamma_1, \ldots , \gamma_m \mbox{ \tiny non-intersecting } \}} \\
\le \sum_{\gamma_1:  l_1 \in R ( \gamma_1) }  e^{ - \beta | \gamma_1 |} \sum_{ \gamma_2 : \min \{ l_i : l_i \notin R (\gamma_1) \} \in R ( \gamma _2 ) }
 e^{ - \beta | \gamma_2|} 
 \ldots  \\
\sum_{ \gamma_m : \min \{ l_i : l_i \notin R ( \gamma_1 ) \cup \ldots \cup R ( \gamma_{ m-1 } ) \} \in R ( \gamma_m)}  e^{ - \beta | \gamma_m |} .
\end{multline}

Note that 
\begin{equation}\label{eq:ub2}
\sum_{ \gamma : l_i \in R( \gamma )} e^{ - \beta | \gamma | } = \sum_{ \gamma : 0 \in R ( \gamma ) } e^{ - \beta | \gamma |} .
\end{equation}

Hence we can upper bound the right hand side of (\ref{eq:lastdisplay}) by 
\begin{multline*}
\sum_{\gamma_1:  l_1 \in R ( \gamma_1) }  e^{ - \beta | \gamma_1 |} \sum_{ \gamma_2 : \min \{ l_i : l_i \notin R (\gamma_1) \} \in R ( \gamma _2 ) }
 e^{ - \beta | \gamma_2|} 
 \ldots \\  \sum_{ \gamma_m : \min \{ l_i : l_i \notin R ( \gamma_1 ) \cup \ldots \cup R ( \gamma_{ m-1 } ) \} \in R ( \gamma_m)}   e^{ - \beta | \gamma_m |} \; \le  \; \left( \sum_{ \gamma : 0 \in R ( \gamma ) } e^{ - \beta | \gamma |} \right)^m  , 
\end{multline*}
where we have applied successively the upper bound (\ref{eq:ub2}) to the right hand side in (\ref{eq:lastdisplay}), starting with $\gamma_m .$ 

To conclude the proof of the lemma, we need an upper bound for  the sum $  \sum_{ \gamma : 0 \in R ( \gamma ) } e^{ - \beta |\gamma |}  .$
Recall that the number of closed contours of length $l $ that contain $0$ is upper bounded by $ 4 l  3^{ l - 2 } .$ Hence,
$$  \sum_{ \gamma : 0 \in R ( \gamma ) } e^{ - \beta |\gamma |}  \le \sum_{ l \geq 4 } 4 l 3^{ l-2} e^{ - \beta l } .$$
Moreover, for $ \beta $ sufficiently large, $  4 l 3^{ l-2} e^{ - \beta l } \le e^{ -  \beta l /2  }.$ Hence
\begin{equation}  \sum_{ \gamma : 0 \in R ( \gamma ) } e^{ - \beta |\gamma |}  \le \frac{ e^{ - 2 \beta }}{ 1 - e^{ - \beta /2 }} \le 1 . 
\end{equation}
We conclude that 
\begin{equation}
\mu_{\beta, \Lambda}^+ \left( x_C = - 1 , x_{\Gamma \setminus C } = + 1 \right) \le  \sum_{ m= 1 }^n    e^{ - 2 \beta n } e^{ - 2 \beta |C|} 
=    n e^{ - 2 \beta n }  e^{ - 2 \beta |C|} \le e^{ - 2 \beta | C|}, 
\end{equation}
for $\beta $ sufficiently large.
Using (\ref{eq:o}), this yields
$$
\P^+_{\beta, \ge , \Lambda} \left( \bigcap_{ i \in \Gamma} \{ X_i = - 1 \} \right)  \le \sum_{ C \subset \Gamma } e^{ - 2 \beta |C|} \ge^{ | \Gamma | - |C| } = \left( e^{ - 2 \beta} + \ge  \right)^{ |\Gamma| }
 .
$$

Letting $\L \to \Z^2 ,$ this concludes the proof of the lemma.
\end{proof}

We are now able to give the proof of Theorem \ref{theo:main}. 

{\bf Proof of Theorem \ref{theo:main}}\\
Let $\Gamma $ be a (self-avoiding) path starting at one of the four $L^1-$neighbors of the origin. We call this path open if $ X_i = - 1 $ for all $i \in \Gamma .$ If $|C_0 ( X) | = \infty ,$ then there exist open paths of all lengths starting at one of the four $L^1-$neighbors of the origin. Write $N(n)$ for the number of such open paths of length $n.$ The number of such possible paths can be bounded from above by $ 4 3^{ n- 1} .$ Thus for any $ n \geq 1,$ using Lemma \ref{prop:1},  
\begin{eqnarray*}
\P_{\beta , \ge}^+ ( | C_0 (X) | = \infty ) & \le  &\P_{\beta , \ge}^+ ( N(n) \geq 1 ) \\
&\le& \E_{\beta , \ge}^+ ( N(n) ) 
= \sum_{ \Gamma : | \Gamma | = n } \P_{\beta , \ge }^+ \left( \bigcap_{i \in \Gamma } \{ X_i = - 1 \} \right)\\
& \le & 4 3^{n-1} 
\left( \ge  + e^{ - 2 \beta  } \right)^n      ,
\end{eqnarray*}
and this converges to $0$ as $n \to \infty ,$
if $ 2 \beta > \ln 3 + e^{- \beta }  $ and $\ge < \frac13 - e^{ - 2 \beta } .$ 

Concerning the proof of item 2., observe that 
\begin{eqnarray*}
\P_{\beta , \ge}^- ( \exists i : | C_i (X) | = \infty ) &\geq&  \mu_\beta^{-} (\exists i : | C_i (X) | = \infty ) \\
&=& \mu_\beta^{-} (\exists \mbox{ an infinite open path } ) = 1 ,
\end{eqnarray*}
by Russo 's classical results (see Proposition 1 of Russo (1979)). 
This concludes the proof.

\section{Proof of Theorem \ref{theo:two}}\label{section3}
For two probability measures $\mu $ and $\nu $ on $A^{\Z^2},$ write $ \mu \preceq \nu $ if there exists a coupling $\bar Q $  having $\mu$ as first marginal and $\nu$ as second marginal,
such that $ \bar Q ( \{ ( x^1 , x^2 ) \in S^2 : x^1 (i) \le x^1 (i) \; \forall i \in \Z^2 \} ) = 1 .$  

Note that 
$$  \P^\ge ( X_i = + 1 | X_j , j \neq i   ) \geq ( 1 - \ge) \lambda_0^+ .$$ 
Now we can apply a standard coupling argument, see for instance Lemma 1.1 of Liggett et al. (1997), to prove that 
$$ \nu_{ ( 1 - \ge) \lambda_0^+} \preceq  \P^\ge \quad  \mbox{   and }  \quad
 \P^\ge \preceq \nu_{ ( 1 - \ge) ( 1 - \lambda_0^-)} .$$
Therefore, 
\begin{multline*}
\P^\ge \left( | C_i ( X)| = \infty \right) \\
= \P^\ge ( \mbox{ there exists an infinite path of $ - 1$ starting from one of the four neighbors of $i $} ) \\
\le \nu_{ ( 1 - \ge) \lambda_0^+} ( \mbox{ there is an infinite path of $ - 1$ starting from one of the neighbors of $i $} ) ,
\end{multline*}
which equals zero by condition (\ref{eq:conditionbetapetit}).

In the same way, under condition (\ref{eq:conditionbetapetit2}),
$$ \P^\ge \left( | C_i ( X)| = \infty \right) \geq \nu_{ ( 1 - \ge) ( 1 - \lambda_0^-)} \left( | C_i ( X)| = \infty \right)  > 0 .$$ 
By Kolmogorov's $0-1$-Law applied to the product measure $\nu_{ ( 1 - \ge) ( 1 - \lambda_0^-)}  ,$  this implies that 
$$ \nu_{ ( 1 - \ge) ( 1 - \lambda_0^-)}\left(\exists i :  | C_i ( X)| = \infty \right) = 1  ,$$
and hence
$$ \P^\ge \left(\exists i :  | C_i ( X)| = \infty \right) = 1 .$$

\section*{Acknowledgments}
We thank two anonymous referees whose remarks helped us to
significantly improve the manuscript. We thank D. Y. Takahashi and
R. Fern\'andez for stimulating discussions and bibliographic
suggestions.  This work is part of USP project MaCLinC, ``Mathematics,
computation, language and the brain", USP/COFECUB project ``Stochastic
systems with interactions of variable range'' and CNPq project
476501/2009-1. It was partially supported by CAPES grant
AUXPE-PAE-598/2011.  A.G. is partially supported by a CNPq fellowship
(grant 305447/2008-4). E.L. has been supported by
ANR-08-BLAN-0220-01. M.C. and E.L. thank NUMEC, University of Sao
Paulo, for hospitality and support.

\vskip30pt
\newpage

Marzio Cassandro

Dipartimento di Fisica

Universit\' a di Roma La Sapienza

P.le A. Moro

00185 Roma, Italy

e-mail: {\tt cassandro@roma1.infn.it}

\bigskip

Antonio Galves

Instituto de Matem\'atica e Estat\'{\i}stica

Universidade de S\~ao Paulo

Caixa Postal 66281

05315-970 S\~ao Paulo, Brasil

e-mail: {\tt galves@usp.br}
\bigskip

Eva L\"ocherbach

CNRS UMR 8088

D\'epartement de Math\'ematiques

Universit\'e de Cergy-Pontoise

95 000 CERGY-PONTOISE,  France

email: {\tt eva.loecherbach@u-cergy.fr}

\end{document}